\documentclass[12pt]{article}
\usepackage{amsmath,amsthm,amsfonts,amssymb,amscd}
\usepackage{graphicx}
\usepackage{euscript}
\def\phi{\varphi}

\newcommand{\be}{\begin{eqnarray}}
\newcommand{\ee}{\end{eqnarray}}

\begin{document}

 \centerline { \bf LARGE DEVIATIONS OF
$U$-EMPIRICAL} \centerline { \bf KOLMOGOROV-SMIRNOV TESTS, AND THEIR
EFFICIENCY }

\vspace{20pt}

 \centerline {  Yakov  Nikitin \footnote{Partially supported by RFBR grant No. 07-01-00159,
and by grant NSh. 638.2008.1.} }

 \vspace{12pt}

 \centerline{ \it St.Petersburg State University}

  \vspace{15pt}

{\footnotesize
 Non-degenerate $U$-empirical
Kolmogorov-Smirnov tests are studied and their large deviation
asymptotics under the null-hypothesis is described. Several examples
of such statistics  used for testing goodness-of-fit and symmetry
are considered. It is shown how to calculate their local Bahadur
efficiency.}

\vspace{15pt}

\section { Introduction.}

Let $X_1,X_2,...$ be i.i.d. observations with continuous
distribution function (df) $F.$ Denote by $F_n$ the usual empirical
df (edf) based on  the first $n$ observations. The Kolmogorov
statistic is the distribution-free statistic defined by
\begin{equation} \label{Kolm} D_n = \sup_{t} | F_n (t) - F(t) |.
\end{equation}

In order to describe the large deviation asymptotics of  statistic
(\ref{Kolm})\ consider for $0<a<1$ the function
\medskip
$$
\begin{array}{ll}
f(a,t)=\begin{cases}
   \displaystyle (a+t)\,\ln\,{\frac {a+t}{t}}+(1-a-t)\,\ln\,{\frac
      {1-a-t}{1-t}}\,, & 0\le t\le  1-a\,,   \\
     &\\
   +\,\infty\,, & 1-a<t\le    1\,,    \end{cases}
\end{array}
$$
and put
   $$ f_0(a)=\inf\limits_{0 < t < 1}\,f(a,t)\,. $$

The following theorem was proved in \cite{Abr}, see also
\cite{Bah1971} and \cite{Niki1995}.

{\bf Theorem 1.1.} \, {\it  For any $ a \in (0,1)$ we have
$$ \label{LD1} \lim_{n \to \infty} n^{-1} \ln {\mathbb P} ( D_n
> a) = - f_0 (a),$$ where the function $f_0(a)$ is continuous on $(0,1),$ and as $ a \to 0$
$$
f_0 (a)=2 a^2 ( 1+o(1))\,.$$ } In particular, this result can be
used to calculate Bahadur efficiency of various modifications of
Kolmogorov-Smirnov tests, see
 \cite{Bah1971} and \cite{Niki1995}.

We are interested in $U$-empirical generalizations of Theorem 1.1
and their applications to  the problem of testing goodness-of-fit
and symmetry. Let $h(x_1,...,x_m)$ be  a real-valued symmetric
kernel of degree $m \geq 1$. Consider the $U$-empirical df (udf)
$$ G_n(t) = {n \choose m}^{-1} \sum_{1\leq i_1 <...<i_m\leq n} I\{
h ( X_{i_1},..., X_{i_m}) < t \}, \quad  t \in R^1.
$$ The properties of such udf's and their use in Statistics were studied in \cite{SER,HJS,
 Jan,SW}.  Also denote  \begin{equation} \label{G} G(t) = {\mathbb P}(h ( X_{1},..., X_{m}) < t
)\end{equation}  and assume that this df is continuous and (for
simplicity) strictly monotonic. Then the $U$-empirical analog of the
Kolmogorov statistic (\ref{Kolm}) has the form
 $$ DU_n = \sup_{t} | G_n(t) - G(t)| $$  and coincides with the Kolmogorov statistic when
  $m=1$ and $h( t ) = t.$

 Various  tests of goodness-of-fit and symmetry can be constructed utilizing the
 Kolmogorov-Smirnov-type statistic
 \begin{equation}
 \label{DDU}
  SU_n = \sup_{t}|G_n (t) - F_n (t)|,  \end{equation}
  with various choices
 of the kernel $h.$   Statistics of this type  frequently appear in constructing statistical
 tests based on characterizations of distributions; then it is  assumed that $G(t) \equiv
 F(t)$, see \cite{Ang,BarHen,BH,niki96}.

The next example illustrates the idea of building such tests.
Consider scale-free testing of exponentiality based on Desu's
characterization \cite{desu}: {\it Let $X_1$ and $X_2$ be
independent non-degenerate and non-negative rv's with common df $F.$
Then the rv's $2 \min (X_1,X_2)$ and $X_1$ have the same
distribution if and only if $F(x) = 1 - \exp(-\lambda x), x\geq 0$
for some $\lambda
> 0.$}

Let us compare the $U$-empirical df
$$ \tilde{G}_n(t) = { n \choose 2}^{-1} \sum_{1\leq j <k \leq n}
I\{ 2\min(X_j,X_k) <t\} $$ with the usual edf $F_n(t)$  by
considering the statistic $$ \label{Desu} DE_n = \sup_t
|\tilde{G}_n(t) - F_n(t)| $$ and assuming its large values to be
critical. The limiting distribution of this statistic can be found
using the results of \cite{Silv}  where weak convergence of $U$-
empirical processes is studied, see also \cite{Denk} and \cite{SER}.
Critical values of $DE_n$ can be calculated via simulation.

 The logarithmic large deviation
asymptotics of $DE_n$ enables  to calculate its Bahadur efficiency
and compare it with other statistics used for testing
exponentiality, see \cite{Bah1971,Niki1995, Serf}. Various examples
of similar statistics will be given  at the end of this paper.

The problem of large deviation asymptotics for $U$-statistics is
studied insufficiently. The abstract large deviation principle
stated in \cite{Eichel} and \cite{SW} is non-effective for
statistical applications. The point is that the rate of decrease for
large deviation probabilities is given implicitly, as a value of
complicated extremal problem in which the Kullback-Leibler
information is minimized on an involved set of probability measures.
The result that could actually be applied to the efficiency
calculations was first obtained for bounded kernels in
\cite{Niki1999}.
 This result is stated as follows.

{\bf Theorem 1.2.} \, {\it  Consider the $U$-statistic of degree $m
\geq 1$
$$
  U_n = {n \choose m}^{-1} \sum_{1\leq i_1 <...<i_m\leq n}\Phi(X_{i_1},...,
X_{i_m}) $$ with centred, bounded, and non-degenerate  real-valued
kernel
 $\Phi$, so that
 $$\label{M} E\Phi(x_1,\dots,X_m)=0, \quad |\Phi(x_1,...,x_m)| \le M,
$$ and  $ \label{sigma}
 \sigma^2 = E\varphi^2 (X_1) > 0
 $
with $ \varphi(y)=E(\Phi(X_1,...,X_m)|X_1=y). $ Then we have
 \begin{equation} \label{limU}
\lim_{n \to \infty} n^{-1} \ln {\mathbb P}( U_n \ge a ):= -
g(a|\Phi) = -\sum_{j=2}^\infty b_j a^j, \end{equation} where the
series converges for sufficiently small $a>0$ and $ b_2 =
(2m^2\sigma^2)^{-1}. $}

The aim of the present paper is to extend Theorem 1.1 (using Theorem
1.2) to the case of statistics $DU_n,$ $SU_n,$ their one-sided
variants, and some other similar statistics.  In the last section we
show how to use the large deviation results for calculating local
Bahadur efficiencies.

\section{Auxiliary arguments and formulation of results}

 To study large deviations of Kolmogorov-type
statistics, it is natural to begin with  a simple one-sided
statistics, see \cite{Abr}, \cite{Bah1971},
$$ DU_n^{+} = \sup_{t} [ G_n(t) - G(t)], \quad DU_n^{-} = \sup_{t} [
G(t) - G_n(t)]. $$ For any $\varepsilon >0$ denote $
P_n(\varepsilon)= {\mathbb P}( DU_n \geq \varepsilon),$ and let
$P_n^{+}(\varepsilon) $ and $P_n^{-}(\varepsilon)$ be the
corresponding probabilities for $ DU_n^{+}$ and $DU_n^{-}.$  It is
evident that
\begin{equation} \label{inii} \max(P_n^{+}(\varepsilon),\
P_n^{-}(\varepsilon)) \leq \ P_n (\varepsilon) \leq \ 2 \max
(P_n^{+}(\varepsilon), P_n^{-}(\varepsilon)).
\end{equation} Therefore, if we prove that for some functions $g^+$
and $g^-$
$$ \lim_{n \to \infty}n^{-1}\ln P_n^{+} (\varepsilon) = - g^+(\varepsilon), \quad
\lim_{n \to \infty}n^{-1}\ln P_n^{-} (\varepsilon) = -
g^-(\varepsilon),
$$
then, by (\ref{inii}), we obtain \begin{equation} \label{plusminus}
\lim_{n \to \infty}n^{-1}\ln P_n(\varepsilon) = - \min(
g^+(\varepsilon), g^-(\varepsilon)). \end{equation}

 The above argument is also valid for the statistic $SU_n$ when we
use one-sided statistics $SU_n^{+}$ and $SU_n^{-}$.  For this
reason, below we consider large deviations of one-sided test
statistics only.

Let  us impose some restrictions on the kernel. We see that the
statistic $DU_n^{+}$ is, in fact, the supremum of the family of
$U$-statistics with the kernels
\begin{equation}
\label{Theta}
 \Theta(x_1,\dots,x_m; t) = I\{ h (x_1,...,x_m) < t \}
- G(t), \quad t \in R^1,
\end{equation} depending on $t.$ The same is true for the
statistic $SU_n^{+},$ where the kernels  also depend on $t$ and have
the form
\begin{equation}
\label{Psi}
 \Psi(x_1,\dots,x_m; t) = I \{ h (x_1,...,
x_m) < t \} - m^{-1} \sum_{i=1}^m I\{x_i < t\}, \quad t \in R^1.
\end{equation}
For any $t$ the kernels $\Theta(\cdot\ ; t)$ and $\Psi(\cdot\ ; t)$
are centred and bounded. Consider their projections
$$ \theta(s; t)=E(\Theta(X_1,...,X_m;t)|X_1=s), \quad  \psi(s; t) =
E(\Psi(X_1,...,X_m;t)|X_1=s),$$ and the variance functions
$$
\sigma_{\theta}^2(t) = E\theta^2(X_1;t), \quad \sigma_{\psi}^2(t) =
E\psi^2(X_1;t).
$$

In the sequel, we will work with general families of kernels
$\Phi(\cdot \ ; t \in {\cal T}),$ that comprise, as particular
cases, the families $\Theta(\cdot \ ; t)$ and $\Psi(\cdot \ ; t)$
defined by (\ref{Theta}) and (\ref{Psi}). For definiteness, we
assume that the parameter set ${\cal {T}}$ is some finite or
infinite interval $[a,b]$ of the real line.

Most results in the literature related to large deviations of
U-statistics deal with U-statistics with specific kernels, see, for
example, \cite{Arc, Eichel, Niki1999, Barin}.  Unlike these results,
we study whole {\it families} of kernels and the corresponding {\it
families} of U-statistics indexed by real parameter $t.$ This is a
distinctive feature of the paper.

We say that the family of $U$-statistics $\{U_n(t), t \in [a,b]\}$
with the kernels $\Phi(\cdot\ ;t) $ and projections $\varphi(\cdot\
;t)$ is {\it non-degenerate, if its variance function
$\sigma^2_{\varphi}(t)= E\varphi^2(X_1; t)$ can vanish only at the
ends of the interval $[a,b]$ and at a finite number of points in the
interior of this interval.}

For example, the one-sided variant $DE_n^+$ of Desu's statistic
(\ref{Desu}) has the variance function equal to $ \frac14 \exp(-t)
(1-\exp(-t)), t \geq 0.$ Another typical variance function  that
appears below in Example 5 is equal to $\frac14 t^2(1-|t|), -1\leq t
\leq 1. $ Both families of $U$-statistics are non-degenerate
according to our definition.

Denote by $\varphi_0^2$ the maximum of the variance function, i.e.,
$\varphi_0^2 =\sup_{t} \sigma_{\varphi}^2(t) > 0$ and introduce the
following subset of the parameter set $:$  $$T=\{ t |\ t\in [a,b],\
\sigma_{\varphi}^2(t)> \frac12 \varphi_0^2\}.$$ On this set the
variance function is separated from zero.  Denote by $T^c $ the
complement of $T$ in $[a,b]$.

Consider the family of $U$-statistics $\{U_n(t), t \in [a,b]\ \}$
corresponding to the family of kernels $\{\Phi(\cdot\ ;t)\}.$  For
$t \in T$ the kernels are non-degenerate, centred, and bounded. Then
the application of Theorem 1.2 yields relation (\ref{limU}) with the
function of two arguments
\begin{equation} \label{gg}g(a,t|\Phi) = \sum_{j=2}^\infty b_j(t)
a^j
\end{equation}
where the coefficients $b_j(t)$ depend on $t$ unlike the function
$g(a|\Phi)$.
 Put \begin{equation}
\label{def} g_T (a|\Phi):= \inf_{ t \in T} g(a,t|\Phi).
\end{equation} It turns out that this function determines the large
deviation asymptotics of the statistics under consideration, and
that the behavior of the function $t \to g(a,t|\Phi)$ on the set
$T^c$ is of no importance. The set $T$ can be somewhat enlarged or
narrowed without essential changes in the results.

As usually in large deviation theory, deriving the large deviation
asymptotics  will consist in obtaining  upper bound and lower bound.
Asymptotic coincidence of the bounds would give the desired result.
We begin with the lower bound.

{\bf Theorem 2.1.} \, {\it Suppose that the non-degenerate family of
kernels $\{\Phi(x_1,...,x_m; t),\\ t\in [a,b] \}$ is bounded and
centred for all $t.$ Then for $a>0$ we have
\begin{equation} \label{low}
 \liminf_{n \to \infty} n^{-1}\ln
{\mathbb P}( \sup_t U_n(t)
>a) \geq - g_{T}(a|\Phi). \end{equation}
The function $a \to g_{T}(a|\Phi)$ is continuous for sufficiently
small $a$ and admits the representation
 \begin{equation}
\label{expan}
 g_{T}(a|\Phi) = \frac{a^2}{2m^2\varphi_{0}^2} +O(a^3), \,\ a
\to 0.\end{equation} }

    From this theorem the lower bound on the large deviation probabilities for
 $U$-empirical Kolmogorov-Smirnov tests follows.  Both one-sided
statistics $DU_n^+$ and $DU_n^-$ satisfy the conditions of Theorem
2.1. They correspond to different families of kernels $
\Theta(\cdot;t)$ and $- \Theta(\cdot;t)$ with  common variance
function $\sigma^2 _{\theta}(t).$ The corresponding functions
$g_{T}(a|\Theta)$ and $g_{T}(a|- \Theta)$ are also different, but
the first terms of their asymptotic expansions as $a \to 0$ are the
same. Now we obtain from (\ref{low})
\begin{equation} \label{DDDU}
 \liminf_{n
\to \infty} n^{-1}\ln{\mathbb P}(DU_n^{\pm}
>a) \geq - g_{T}(a|\pm \Theta),
\end{equation} and  similarly  \begin{equation} \label{SSSU} \liminf_{n
\to \infty} n^{-1}\ln {\mathbb P}( SU_n^{\pm}
>a) \geq - g_{T}(a|\pm \Psi).
\end{equation}

For the purpose of obtaining the upper bound, we assume that the
following condition of {\it monotonicity in parameter} on the family
$\{U_n(t), t\in[a,b]\}$ is satisfied. Suppose there exists a
sequence of partitions of the interval $[a,b]$ into $N$ parts:
$a=t_0<t_1< \dots < t_N=b,$ such that the nodes of the partition do
not coincide with the zeros of the variance function and that for
any $k = 0,\dots, N-1$
\begin{equation} \label{MC} \sup_{t_k \leq t < t_{k+1}}U_n(t) \leq
U_n(t_{k+1}) + \Delta_n (N),\end{equation} where the sequence of
rv's $\Delta_n (N)$  {\it decreases fast}, when $n$ and $N$ grow.
More precisely, it means that there exists a sequence
 $\{\tau_N\}, \tau_N \to 0\ \,\mbox{as}\,\, N \to \infty,$ such
 that
 \begin{equation} \label{fast} \lim_{N
\to \infty} \lim_{n\to \infty} n^{-1}\ln {\mathbb P}(\Delta_n
(N)>\tau_N) = -\infty. \end{equation}

Under this condition the upper bound result takes the following
form.

{\bf Theorem 2.2.} \, {\it Suppose that the conditions of Theorem
2.1 are valid and that the condition of monotonicity in parameter of
the family $\{U_n(t)\}$ holds. Then we have
\begin{equation} \label{4} \limsup_{n \to \infty}\ n^{-1}\ln {\mathbb
P}( \sup_t U_n(t)
>a) \leq - g_{T}(a|\Phi). \end{equation}}

    Combining Theorems 2.1 and 2.2, we arrive at the following result.

{\bf Theorem 2.3.} \, {\it Under the conditions of Theorem 2.2 we
have
$$ \lim_{n \to \infty}
n^{-1}\ln {\mathbb P}( \sup_t U_n(t)
>a) = - g_{T}(a|\Phi). $$}

It is shown below that both statistics $DU_n^{\pm}$ and $SU_n^{\pm}$
satisfy the monotonicity condition. If the corresponding families of
kernels (\ref{Theta}) and (\ref{Psi}) are non-degenerate and
centred, then Theorem 2.1 is  applicable.  Hence (\ref{inii}) and
(\ref{plusminus}) imply the following result for the two-sided
statistics.

{\bf Theorem 2.4.} \, {\it  If the family of kernels {\rm
(\ref{Theta})} is non-degenerate, then
$$
\lim_{n \to \infty} n^{-1}\ln {\mathbb P}( DU_n
>a) =  v_D(a)= - \min( g_{T}(a|\Theta), g_{T}(a|-\Theta)). \\
$$
If the family of kernels {\rm (\ref{Psi})} is non-degenerate and
centred, then
$$
 \lim_{n \to \infty} n^{-1}\ln {\mathbb P}( SU_n >a) = v_S
(a)= - \min(g_{T}(a|\Psi), g_{T}(a|-\Psi)),
$$
The functions $v_D$ and $v_S$  are continuous and satisfy the
asymptotic relations
$$ v_D (a) \sim - a^2/2m^2\theta_0^2,\quad v_S (a) \sim
-a^2/2m^2\psi_0^2,\quad  a \to 0.
$$}

\section{Lower bound.}

In order to prove the lower bound we need the theorem on implicit
analytic operators. Consider three Banach spaces $E_1$, $E_2$, and
$E_3$. Denote by ${\cal D}_r(x_0,E)$ the ball in the space $E$ of
radius $r$ with center at $x_0$. We are interested in finding
solutions $x=x(y)$ of the operator equation
\begin{equation} \label{eqq} F(x,y)=0 \end{equation} (here $y$ plays a part of
a parameter)  under the condition \begin{equation} \label{boundary}
x(y_0)=x_0. \end{equation}  We assume that the operator $F(x,y)$ is
analytic and the condition $F(x_0,y_0)=0$ holds.   For main
definitions and facts of the theory of analytic operators in Banach
spaces, we refer to \cite[ \S 22]{VT} and \cite[\S 32]{Kras}. In
particular, the operator $F(x,y)$ is called analytic in some domain
if in the neighbourhood of any point it can be represented as an
uniformly convergent Taylor operator series \cite{VT}.

{\bf Theorem 3.1.} \, {\it see \cite[Theorem 22.2 ]{VT}. Suppose
that the operator $F(x,y)$ is analytic in ${\cal D}_r(x_0,E_1)\times
{\cal D}_\rho(y_0,E_2)$ with values in $E_3$. Let the operator
$B=F^\prime_x(x_0,y_0)$ (the derivative is understood in the
Fr\'{e}chet sense) has a bounded inverse operator. Then there exist
positive numbers $r_1$ and $\rho_1$ such that the equation {\rm
(\ref{eqq})} has an unique solution $x= f(y)$ in the ball ${\cal
D}_{r_1}(x_0,E_1)$. This solution is defined and analytic in the
ball ${\cal D}_{\rho_1}(y_0,E_2),$ and satisfies the condition {\rm
(\ref{boundary})}.}

We will prove Theorem 3.1  using the arguments of \cite{Niki1999} in
conjunction with some auxiliary results. The main idea  is to
construct a majorant series for solutions of nonlinear equations
that would guarantee the uniform convergence with respect to the
parameter $t.$

{\it Proof of Theorem 3.1.} It is clear that for any statistic of
the form $\sup_t U_n(t),$ with kernels $\Phi(\cdot\ ;t)$ and
projections $\varphi(\cdot\ ;t),$ the following holds true:
$$ {\mathbb P}(\ \sup_{t} U_n(t)
>a) \geq \sup_{t\in T } {\mathbb P} (U_n(t) >a).
$$
Consequently, $$ \liminf_{n \to \infty} n^{-1} \ln {\mathbb P} (
\sup_t U_n (t)
> a) \geq -\inf_{t\in T} g(a,t|\Phi) = - g_T (a|\Phi).
$$
Hence in order to obtain the lower bound we must analyze the
function $ g(a,t|\Phi)$  in (\ref{gg}) over the set $T.$ Recall that
for $t\in T$ we have $\sigma_{\varphi}^2(t)> \frac12 \varphi_0^2
>0.$

For any Borel set  $S \subset R^1$ consider the Banach space of
measurable bounded functions $B(S)$ with the norm
$$
\parallel x \parallel = {\sup}_{s\in S} |x(s)| \ .
$$

We may assume that the initial  observations $X_i,i =1,\dots,n,$ are
uniformly distributed on $I= [0,1].$ Otherwise we can consider the
sample $U_i= F(X_i),i=1,\dots,n,$ and replace the initial kernel
$h(X_1,...,X_m)$ by $h(F^{-1}(U_1),\dots,F^{-1}(U_m)).$ The families
of kernels $\Phi(\cdot\ ;t), \Psi(\cdot\ ;t), \, \Theta(\cdot\ ;t) $
and the corresponding families of $U$-statistics depending on $U_i$
remain centred, bounded, and non-degenerate.

For simplicity we consider only the kernels of degree 2. With slight
changes in the proof, the results  remain valid for kernels of any
finite degree.

For any $t \in T$ we use the variant of Sanov's theorem for large
deviations of $U$-empirical measures from \cite{Eichel}, see also
\cite{SW}. So, we reduce the problem of large deviations to the
problem of the minimization of entropy under suitable normalization
conditions. It follows from \cite{Niki1999} that the function of
interest $g(a,t|\Phi)$ is the solution of the extremal problem
\begin{equation} \label{Extr} g(a,t|\Phi)= \inf\{ \int_0^1 (1+a
x(s;t)) \ln (1+ a x(s;t))ds:
  x \in B ( I\times T), \int_0^1x(s;t)ds =0\}\end{equation} under the normalization condition
 \begin{equation} \label{norma}
 2\int_0^1\varphi(s_1;t)x(s_1;t)ds_1 + \int_0^1\int_0^1
\Phi(s_1,s_2;t) x(s_1;t)x(s_2;t)ds_1ds_2=1. \end{equation}

The Euler-Lagrange equation for the extremal $x(s;t)$ from
$B(I\times T)$ takes the form, see \cite{Niki1999} :
\begin{equation}
\begin{array}{rl}\label{inteq}
 (1+a x(s_1;t))\int_0^1 \exp\left\{\lambda(t) a \int_0^1
\Phi(s_1,s_2;t) (1+a x(s_2;t))ds_2 \right\}ds_1-\\   -\exp
\left\{\lambda(t) a\int_0^1 \Phi(s_1,s_2;t) (1+a x(s_2;t))ds_2
\right\} =0,
\end{array}
 \end{equation}
under the same normalization condition (\ref{norma}). To simplify
the notations denote
$$ \rho(x|s_1;t) = \int_0^1 \Phi
(s_1,s_2;t)x(s_2;t)ds_2. $$ Then (\ref{inteq}) can be written as
\begin{equation}
\begin{array}{rl} \label{XX} (1+a x(s_1;t))\int_0^1
\exp\left(\lambda(t) a \varphi(s_1;t) +\lambda(t) a^2
\rho(x|s_1;t)\right)ds_1- \\  -\exp\left(\lambda(t) a \varphi(s_1;t)
+\lambda(t) a^2 \rho(x|s_1;t)\right)=0,
\end{array}
\end{equation}
and the condition (\ref{norma}) becomes
 \begin{equation}
\label{norm2}
 \int_0^1 (2\varphi(s_1;t) + a\rho(x|s_1;t))x(s_1;t)ds_1 =1.
 \end{equation}

Expanding the exponents on the left-hand side of (\ref{XX})  into a
series and integrating, we get
$$
\begin{array}{rl}
\vspace{8pt} &\sum_{k+j \geq 0}  \lambda^{k+j}(t) a^{k+2j} (k!
j!)^{-1}\int_0^1 \varphi^k (s_1;t)\rho(x|s_1;t)ds_1 - \\
\vspace{8pt}
 &\sum_{k+j \geq 0}  \lambda^{k+j}(t) a^{k+2j}(k! j!)^{-1}
 \varphi^k (s_1;t)\rho(x|s_1;t) + \\
 \vspace{8pt}
 &\sum_{k+j \geq 0}  \lambda^{k+j}(t) a^{k+2j+1}(k!
j!)^{-1}x(s_1;t)
 \int_0^1 \varphi^k (s_1;t)\rho(x|s_1;t)ds_1 =0.
\end{array}
$$

To simplify this equation we extract two first terms corresponding
to the indices $k=0,j=0 $ and $k=1,j=0,$ so that the remainder of
the sum is taken over the set of indices $N(k,j) = \{(k,j) \neq
(0,0),(1,0)\}.$  After some algebra we obtain the equation
\begin{equation}
\begin{array}{rl}\label{16} x(s_1;t) -\lambda(t)\varphi(s_1;t)  +
\sum_{N(k,j)} \lambda^{k+j}(t) a^{k+2j-1}(k! j!)^{-1} \times \\
 \times \left[
 \varphi^k (s_1;t)\rho(x|s_1;t)
- (1+a x(s_1;t))
 \int_0^1 \varphi^k (s_1;t)\rho(x|s_1;t)ds_1 \right]=0.
 \end{array}
\end{equation}

Note that for $a=0$ the ``principal" \ solution of this equation
satisfying the normalization condition is $ x_0(t,s) = \lambda_0(t)
\varphi(s;t),\mbox{where}\,\lambda_0(t) =
(2\sigma_{\varphi}^2(t))^{-1}.$ Our aim is to build the
``perturbation" of this solution for $a
>0.$

We introduce a new small functional parameter $\nu(t) = \lambda(t)
-\lambda_0(t)$ and  a new unknown function $y(t,s) = x(t,s) -
\lambda_0 (t) \varphi(s;t),$ and substitute them into
 equation (\ref{16}). We have \begin{equation}
\begin{array}{rl}
\label{117} \vspace{8pt} &y(s_1;t)- \nu(t) \varphi(s_1;t)+\\
\vspace{8pt}
 &\sum_{N(k,j)} (\nu(t)
+\lambda_0(t))^{k+j} a^{k+2j-1}(k! j!)^{-1}[
 \varphi^k (s_1;t)\rho^j(y+\lambda_0(t) \varphi(s_1;t)|s_1;t)-\\
\vspace{8pt}
  & \int_0^1 \varphi^k (s_1;t)\rho^j(y+\lambda_0(t) \varphi(s_1;t)|s_1;t)ds_1]-
\sum_{N(k,j)} (\nu(t) +\lambda_0(t))^{k+j} a^{k+2j}(k!
j!)^{-1}\times\\
 &(y(s_2;t)+ \lambda_0(t)\varphi(s_2;t))\int_0^1\varphi^k (s_1;t)
 \rho^j(y(s_1,t)+\lambda_0(t)\varphi(s_1;t)|s_1,t)ds_1=0.
\end{array} \end{equation}

Due to the inequality
$$ \lambda_0(t) =(2\sigma_{\varphi}^2(t))^{-1} <
\varphi_{0}^{-2}, \quad t\in T ,$$ and the boundedness of the
kernel, the series on the left-hand side are convergent series of
$k,j$-linear operators (see \cite{VT} ) in $y$ and $\nu$ with
bounded coefficients. Therefore the left-hand side of equation
(\ref{117}) is the analytic operator
$$A(y,\nu, a): B(I\times T)\times B(T) \times R^1
\to B(I\times T).$$

The Fr\'{e}chet derivative  $ A_y(y,0,0) $ at the point $y=0$ is the
unit operator and hence is bounded. Then Theorem 3.1 guarantees the
existence of a solution of the form
\begin{equation} \label{sol}
 y(s_1;t) = \sum_{k+j\geq 1}
c_{kj}(s_1;t) \nu^k(t) a^j, \end{equation} where the series is
absolutely convergent in the space $ B(I \times T)$ for sufficiently
small $||\nu||$ and $a.$  This means the convergence of the power
series with numerical coefficients $
 \sum_{k+j\geq 1}
||c_{kj}||\cdot ||\nu||^k a^j $ for sufficiently small $||\nu||$ and
$a.$ Note that in the proof of Theorem 3.1 in \cite{VT} the majorant
series for the solution was built explicitly.

Now we substitute solution (\ref{sol}) into equation (\ref{117}).
Equating the coefficients at the same powers of $\nu$ and $a$, we
obtain the expressions for $c_{kj}$. For example,
$$
 c_{10}(s;t) = \varphi(s;t),\,\, c_{01}(s;t) =\frac{1}{2} \lambda_0^2(t)
(\varphi(s;t) - 3\sigma_{\varphi}^2(t))+ \lambda_0(t) \rho(s;t),$$
and so on. Returning to the function $x,$ we have
 \begin{equation} \label{xynu}
 x(s;t)
= \lambda_0(t) \varphi(s;t) +  \sum_{k+j\geq 1} c_{kj}(s;t) \nu^k
(t) a^j .\end{equation} Substituting this solution into
normalization condition (\ref{norm2}), we obtain \begin{equation}
\begin{array}{rl}
\label{NNN} \vspace{8pt} &2 \sum_{k+j\geq 1}\int_0^1
c_{kj}(s;t)\varphi(s;t)ds\ \nu^k(t)a^j
\\ \vspace{8pt} &+\sum_{k+j\geq 1} \sum_{i+l\geq 1} \int_0^1
\int_0^1 \Phi(s_1,s_2;t)c_{kj}(s_1;t)c_{il}(s_2;t)ds_1ds_2\ a^{i+j}
\nu^ {k+l}(t)=1. \vspace{8pt}
\end{array}
 \end{equation}
As $c_{10}(s;t) = \varphi(s;t)$, the coefficient at $\nu(t)$ is
equal to $2\sigma^2(t)$ and is positive on $T.$ Dividing by
$2\sigma^2(t),$ we have the equation
 \begin{equation}
 \nu(t) =\sum_{k\geq 1,l \geq 2} \gamma_{kl}(t) a^k\nu^l(t),
\end{equation}
 where, as seen from
(\ref{NNN}), the series with coefficients $\gamma_{kl}(t)$ converges
absolutely in some ball of the space $B(T).$

Applying again  Theorem 3.1 to equation (\ref{NNN}), we obtain the
representation \begin{equation} \label{nuser} \nu(t) = \sum_{p\geq
1} e_p (t) a^p, \end{equation} where the series converges absolutely
in $B(T)$ for sufficiently small $a
> 0.$ Substituting (\ref{nuser}) into (\ref{xynu}), we get
again the convergent series. Returning to the extremal $x(s,t) =
\lambda_0(t) x_0(s;t) + y(t),$ we substitute the new series for $x$
into (\ref{Extr}). Integrating term-wise and using the convergence
of the series for the solution, we obtain the expression for
$g(a,t)$ of the form  \begin{equation} \label{ser} g(a,t| \Phi)=
 \sum_{i=2}^{\infty} c_i(t) a^i,\, t \in T,
 \end{equation}
where $ c_2(t) = (8\sigma_{\varphi}^2 (t))^{-1}$, and the series is
convergent for sufficiently small positive $a$, so that  the
majorant series $ \sum_{i=1}^{\infty} ||c_i|| a^i$ is also
convergent.

Now we can prove the continuity in $a$ and other properties of the
function $g_{T}(a| \Phi)$ listed in Theorem 2.1. Indeed, for any
$a_1
> a_2$ from the interval of convergence of the majorant series we have
 $$
\begin{array}{ll}
\vspace{8pt}
 |g_{T}(a_1|\Phi) - g_{T}(a_2|\Phi)| &\leq
\sup_{t \in T}|\sum_{k=2}^{\infty} c_k(t)(a_1^k -a_2^k)| \\
\vspace{8pt}
 &\leq |a_1
-a_2| \sum_{k=2}^{\infty}k \ ||c_k||a_1^{k-1}.
\end{array}
$$  The series  $\sum_{k=2}^{\infty}k \ ||c_k||a_1^{k-1}$ has the same radius of convergence
as the majorant series $ \sum_{i=1}^{\infty} ||c_i|| a^i$, hence the
sum $\sum_{k=2}^{\infty} k \ ||c_k|| a_1^{k-1}$ is bounded, and the
continuity of $g_{T}(a|\Phi)$ follows.

Now let us estimate the difference $|g_{T}(a|\Phi) - a^2/(8
\varphi_{0}^2)|$ for small $a.$ Note that $$ g_{T}(a|\Phi) =
\inf_{t\in T}g(t,a|\Phi) \leq a^2/8 \varphi_0^2 + \sup_{t\in
T}|\sum_{k\geq 3} c_k(t)a^k|.$$ On the other hand,
$$ g_{T}(a|\Phi) = \inf_{t\in T}g(t,a|\Phi) \geq
 a^2/8\varphi_0^2  - \sup_{t\in T} \sum_{k\geq 3}(-c_k(t))a^k. $$
 Hence we obtain
 $$
\begin{array}{ll}
|g_T(a|\Phi) - a^2/(8\varphi_{0}^2)| \leq \sum_{k=3}^{\infty}||c_k||
a^{k} = a^3 \sum_{k=3}^{\infty} ||c_k|| a^{k-3} = O(a^3),\, a\to 0.
\end{array}
$$

For kernels of degree $m > 2$ in (\ref{expan}), the term
$2m^2\varphi_0^2$ appears instead of $8\varphi_0^2$. $ \hfill
\square $

\section{Upper bound.}

In this section, we bound the large deviation probabilities from
above. For this we apply the exponential inequality for
non-degenerate $U$-statistics from
 \cite[Theorem 2.]{Arco} For simplicity, we  give here a slightly weaker version of
it.

{\bf Lemma 4.1} {\it Under the conditions of Theorem 1.2 with $
\sigma^2
>0$ and $z>0$ we have $$ \label{sArc} {\mathbb P}( |U_n|
> z ) \leq 4\exp\left( - \frac{ n z^2}{2m^2 \sigma^2 +
Lz}\right), $$ where $ L:= \left(2^{m+3} m^m + (2/3)m^{-1}\right)M.
$}

\newpage

{\it Proof of Theorem 2.2.} Consider the partition of the parametric
set $a=t_0<t_1<\dots<t_N=b$ from the monotonicity condition with the
nodes different from the zeros of the variance function
$\sigma^2_{\varphi}(t).$ On any interval of the form $[t_k,
t_{k+1}),\ k=0,\dots, N-1, $ we have  \begin{equation} \label{inin}
 \sup_{t_k \leq t < t_{k+1}} U_n(t) \leq U_n(t_{k+1}) +
\Delta_n(N).  \end{equation} Next,  using (\ref{inin}) for $\tau_N
>0$ from condition (\ref{MC}),
$$
\begin{array}{ll}
\vspace{8pt}
  {\mathbb P}(\sup_{t} U_n(t) > a) &\leq \sum_{k=0}^{N-1}{\mathbb P} (\sup_{t_k \leq t < t_{k+1}} U_n(t)
  \geq a) \leq \\ &\leq
 \sum_{k=0}^{N-2}{\mathbb P} (U_n(t_{k+1})  \geq a - \tau_N) +  N{\mathbb P}( \Delta_n(N) >\tau_N)
 = \Gamma_{1,N}+ \Gamma_{2,N}.
\end{array}
$$
Let us divide the sum $\Gamma_{1,N}$ into two parts: the first sum
includes the indices $k$ for which $t_{k+1} \in T$, while the second
sum includes all remaining indices. For $k=0,\dots, N-2$ the rv $U_n
(t_{k+1})$ is a $U$-statistic with centred and bounded kernel
$\Phi.$  By Theorem 1.2 we have for the summands of the first sum
 $$
 {\mathbb P} ( U_n (t_{k+1}) > a- \tau_N) = \exp(-n g(a- \tau_N,
t_{k+1})+ o(n)) \leq  \exp(-n g_T(a- \tau_N) + o(n)).
 $$
For the summands of the second sum,  by Lemma 4.1 for $t_{k+1}
\notin T$
$$ {\mathbb P} (U_n (t_{k+1}) > a- \tau_N)
 \leq 4\exp\left( - \frac{ n (a- \tau_N)^2}{m^2
\varphi_0^2 + L(a- \tau_N)}\right).
$$
Therefore \begin{equation} \label{esti} \Gamma_{1,N} \leq N\exp(-n
g_T(a- \tau_N|\Phi) + o(n))+ 4N\exp\left( - \frac{ n (a-
\tau_N)^2}{m^2 \varphi_0^2 + L(a- \tau_N)}\right). \end{equation}

 Thanks to (\ref{fast}) the term $\Gamma_{2,N}$  decreases
faster than $\Gamma_{1,N},$ and  can be neglected. Taking the
logarithms in the inequality (\ref{esti}), dividing by $n$ and
passing to the limit as $n \to \infty,$ we obtain $$ \limsup_{n \to
\infty}n^{-1} \ln {\mathbb P} (\sup_t U_n (t) > a) \leq - \min
\left( g_T(a- \tau_N|\Phi), \frac{ (a- \tau_N)^2}{m^2 \varphi_0^2 +
L(a- \tau_N)}\right).
$$
By continuity of the function $g_T(\cdot \ | \Phi)$ as $N \to
\infty$
$$ \limsup_{n \to \infty}n^{-1} \ln {\mathbb P} (\sup_t U_n
(t) > a) \leq - \min \left( g_T(a|\Phi), a^2 (m^2 \varphi_0^2 + L a
)^{-1}\right).
$$
But for small $a,$ by (\ref{expan}), $ g_T(a|\Phi) < a^2 (m^2
\varphi_0^2 + L a )^{-1}.$ Hence we obtain the  required inequality
$$ \limsup_{n \to \infty}n^{-1} \ln {\mathbb P} (\sup_t U_n
(t) > a) \leq -  g_T(a|\Phi).\qquad \qquad \hfill \square$$

Theorem 2.3 follows immediately from Theorems 2.1 and 2.2.

\section{Kolmogorov-Smirnov-type statistics}

In this section we prove that the Kolmogorov-Smirnov-type statistics
$DU_n^{\pm}$ and $SU_n^{\pm}$ satisfy the monotonicity condition, so
that Theorem 2.3 is applicable to them.

First, consider the statistic $DU_n^+$ and assume that the family of
kernels $\Theta (\cdot\ ;t)$ is non-degenerate and $\theta_0^2>0.$
Let $N$ be a large number such that $a- N^{-1}
>0.$ We divide the parametric set into $ N$ parts with nodes
$t_k = G^{-1} (k/N), k=0,...,N,$ where $G$ is from (\ref{G}). If
some interior node coincides with the zero of the variance function,
we make a shift of order $O(N^{-2}).$ On any interval $[t_k,
t_{k+1}), k = 0,...,N-1 $ we have
$$\sup_{t_k \leq t < t_{k+1}} (G_n(t) - G(t)) \leq G_n (t_{k+1}) - G(t_{k+1}) +O(N^{-1}).
$$
Hence the monotonicity condition (\ref{MC}) holds if we take
$\Delta_n(N)= O(N^{-1}),$ and $\tau_N$ to be any sequence tending to
zero slower than $N^{-1}.$  The same procedure is applied to the
statistic $DU_n^-$.

The arguments for $SU_n^{+}$ are similar. Take the nodes $t_{k,N} =
F^{-1}(\frac kN), k=0,...,N, $ and shift them, if necessary, as
above. Hence for $k=0,\dots,N-1$ we have
$$
\sup_{t_k \leq t < t_{k+1}}[G_n(t)- F_n(t)] \leq G_n(t_{k+1}) -
F_n(t_{k+1}) + F_n(t_{k+1})- F_n(t_{k}).
$$
The part of the quantity $\Delta_n(N)$ in the monotonicity condition
plays the rv
$$ F_n(t_{k+1})- F_n(t_{k})= n^{-1}\sum_{j=1}^n I\{t_k \leq
X_j < t_{k+1} \}. $$ Obviously, the sum on the right-hand side has
the binomial distribution with parameters $n$ and $ p= 1/N.$ The
next lemma is proved in \cite{niki96}.

 {\bf Lemma 5.1} {\it Let $Bin(n,1/N)$ be the rv having binomial distribution with
parameters $n$ and $1/N.$ Then for any $\tau \in (0,1) $ the
following inequality holds
  $$ \label{ineq} {\mathbb P}( Bin(n,1/N)
> n\tau) \leq 4^{n} \exp( -n\tau \ln N ). $$}

\smallskip
\noindent  We apply this Lemma with $ \tau= \tau_N= (\ln N)^{-1/2}.
$ For sufficiently large $N$ $$ \label{vazhno} {\mathbb P}(
Bin(n,1/N)
> n / \sqrt{\ln N} ) \leq  \exp(-\frac12 n \sqrt{ \ln N}), $$
and hence  \begin{equation} \label{Chern2} {\mathbb P}(F_n(t_{k+1})-
F_n(t_{k})
>\tau_N)\leq  {\mathbb P}(Bin(n, 1/N) > n\tau_N) \leq  \exp( -
\frac{1}{2} n\sqrt{\ln N}). \end{equation} From this the
monotonicity condition follows. The arguments for $SU_n^-$ are
similar.

We see that Theorem 2.2 is applicable to statistics $DU_n^{\pm}$ and
$SU_n^{\pm}$ if corresponding families of kernels are non-degenerate
and centred. In this case
$$\begin{array}{ll}\vspace{8pt} &\limsup_{n \to
\infty}n^{-1}\ln {\mathbb P}(DU_n^{\pm} > a) \leq - g_{T}(a|\pm\Theta),\\
&\limsup_{n \to \infty}n^{-1}\ln {\mathbb P}(SU_n^{\pm} > a) \leq -
g_{T}(a|\pm \Psi).
\end{array}
$$
Together with (\ref{DDDU}) and (\ref{SSSU}) this implies
$$\begin{array}{ll}\vspace{8pt}
&\lim_{n \to
\infty}n^{-1}\ln {\mathbb P}(DU_n^{\pm} > a) = - g_{T}(a|\pm\Theta),\\
&\lim_{n \to \infty}n^{-1}\ln {\mathbb P}(SU_n^{\pm} > a) = -
g_{T}(a|\pm \Psi).
\end{array}
$$
From these relations, as explained at the end of section 2, Theorem
2.4 follows.

In order to illustrate the result on large deviations of the
statistic $DU_n,$ assume for simplicity that the initial continuous
d.f. $F$ is defined on some finite or infinite interval $[a,b]$ and
is strictly monotonic there. Consider the kernel $h(x,y) =
\max(x,y).$ Then $E I\{{h(X, Y)< t}\} = F^2 (t).$ It follows that
the kernels corresponding to one-sided statistics $DU_n^{(max)}$
have the form $ \Theta(x,y; t)= I\{\max(x,y) < t\} - F^2 (t), t \in
R^1. $ The projections of these kernels are
$$
\theta(x; t) = {\mathbb P}(\max(x,Y)<t) - F^2(t) = I\{x<t\} \ F(t) -
F^2 (t),
$$
and the common variance function is
$$ \sigma^2_{\theta} (t) = E \theta^2(X; t) =
F^3(t)(1-F(t)). $$ The maximum of this function is $27/256.$ It now
follows that for some continuous function $v_0$
$$
\lim_{n\to \infty}n^{-1} \ln {\mathbb P}( DU_n^{\left(\max\right)}
\geq a) = v_0(a) = - \frac{32}{27} a^2 +O(a^3),\, \, a \to 0.
$$

\section{Statistical applications}

In this section we apply our general theorems proved above to
particular $U$-empirical tests of Kolmogorov-Smirnov type. At the
same time, we fill small gaps in the proofs of \cite{BH, niki96} and
\cite{Abb}, where the incorrect paper \cite{das} was used.

1. {\it Test of exponentiality based on Desu's characterization.}
Let us return to the statistic $DE_n$  given in the Introduction.
This statistic is scale-free, so we can assume that the observations
have standard exponential distribution. The kernel of the family of
$U$-statistics $DE_n^+$ takes the equivalent form
$$ \Psi(x,y; t) =  \frac12 ( I\{x >t\} + I\{y>t\})- I\{ \min(x,y) >
t/2 \},\quad x,y,t \geq 0.
$$
  The projection is given by
$$
\begin{array}{ll}
&\psi(y; t)= \frac12 ( {\mathbb P}\{X_1>t\} + I\{y>t\}]  - E[I\{
\min(X_1,y)> t/2 \}]
\\ &=\frac12 \exp(-t) +\frac12 I\{y>t\} -  \exp(-t/2) \ I\{ y > t/2
\},
\end{array}
$$
and hence the variance function $ \sigma_{\psi}^2 (t)
 =\frac14 \exp(-t) (1-\exp(-t)), t \geq 0.$ So, the family of kernels
 is non-degenerate, and we can apply Theorem 2.4.
 It is seen that $\psi_0^2 = \sup_{t\geq 0}\sigma_{\psi}^2 (t)=
\frac{1}{16},$ and this determines the first term of the large
deviation asymptotics. So, there exists a continuous function $v_1$
such that \begin{equation} \lim_{n \to \infty} n^{-1}\ln {\mathbb
P}( DE_n
> a) =  v_1(a) = -2a^2 + O(a^3), \quad a \to 0. \end{equation}

\medskip

 2. Another test of exponentiality is based on the simplified
``lack of memory" property, see \cite{Ang} and \cite{niki96}.
Consider the statistic  $$ AN_n^+ = \sup_{x\ge 0}[\bar F_n(2x) -\bar
F_n^2(x)].
$$
Statistics $AN_n^-$ and $AN_n$ are defined analogously. Large values
of these statistics are statistically significant.

It was shown in \cite{niki96} that under the hypothesis of
exponentiality the statistic $AN_n^+$ admits the representation
$$AN_n^+ = \sup_{0<t<1} {n \choose 2}^{-1} \sum_{1\leq i < j\leq n}
\Psi(U_i,U_j; t) + O(n^{-1}),
$$ where $U_1,...,U_n$ are uniformly distributed on [0,1] rv's. The family of kernels
is then given by
\begin{multline*}
\Psi(x_1,x_2; t)=I\{x_1<t\}+I\{x_2<t\} - I\{x_1<t\}I\{x_2<t\}
-\\(1/2)(I\{x_1<2t-t^2\} + I\{x_2<2t-t^2\}).
\end{multline*}
The corresponding family of $U$-statistics satisfies the
monotonicity condition \cite{niki96}. Simple calculations show that
$$ \psi(z; t)=E\left(\Psi(X_1,X_2; t)|X_1=z \right)=(1-t)I\{z<t\} + t^2/2 -
(1/2)I\{z<2t-t^2\}$$ and hence the variance function is  equal to
$$ \sigma^2_{\psi}(t) =(1/4)t(1-t)^2(2-t)>0,\quad  0<t<1.$$
We can apply Theorem 2.4.

The maximum of the variance function is attained at $t= 1 -
1/\sqrt{2}$ and is equal to $1/16.$ Therefore for some continuous
function $v_2$
$$ \lim_{n \to \infty} n^{-1} \ln {\mathbb P}(AN_n^+ \geq a) = v_2(a) = -2a^2
+O(a^3),\, a \to 0.
$$
The same asymptotics is valid for statistics $AN_n^-$ and $AN_n.$

3. One more characterization of the exponential law belongs to Puri
and Rubin \cite{Puri}: {\it Let $X$ and $Y$ be independent rv's with
common absolutely continuous df $F$ on $R^{+}.$  $F$ is exponential
if and only if  $|X - Y|$ has the same distribution as $X.$}

We construct a one-sided $U$-empirical version of the
Kolmogorov-Smirnov test by introducing the statistic
$$
PR_n^+ = \sup_{t\geq 0} \left({n \choose 2}^{-1} \sum_{1\leq i <
j\leq n} I\{ |X_i - X_j| <t\} - n^{-1}\sum_{k=1}^n I\{ X_k
<t\}\right).
$$
Statistics $ PR_n^- $ and $ PR_n $ are defined analogously. Thus we
obtain the family of $U$-statistics with the kernels
$$\Psi(x_1,x_2; t)= I\{|x_1 -x_2|<t\} - \frac12( I\{x_1 <t\} +
I\{x_2 <t\}),\ t \geq 0. $$  Hence
\begin{multline*}
$$ \psi(s; t) = E\left(I\{|X_1 -s|<t\} -
\frac12( I\{X_1 <t\} + I\{s <t\})\right)=\\= \left( e^{-s +t} -
\frac12 \right)  I\{s\geq t\} - e^{-t}(e^{-s}
  - \frac12),
\end{multline*}
and the variance function is equal to
$$
\sigma_{\psi}^2 (t) = \frac{1}{12}e^{-t}(1+e^{-t}-2e^{-2t}), \,
t\geq 0.
$$

The maximum of this function is attained for $ e^{-t} =
\frac{\sqrt{7}+1}{6} $ and is equal to $ \frac{10+7\sqrt{7}}{648}.$
Hence there exists a continuous function $v_3$ \ such that as $a \to
0$
$$ \lim_{n \to \infty} n^{-1}\ln {\mathbb P}( PR_n
> a) =  v_3(a) = -\frac{7\sqrt{7}-10}{3} a^2 + O(a^3)\approx  -
2.840 a^2 + O(a^3).
 $$

The next two examples are related to the problem of testing
symmetry.

4. Let $X_1, X_2, \dots, X_n$ be a random sample from a continuous
df $F_0$ . We wish to test the hypothesis of symmetry  about zero
$$H_0: \Delta F_0(x):= F_0(x)+F_0(-x)-1=0 \quad \forall x \in
{\mathbb R}^1.$$  Consider the centered (in time) Kolmogorov-Smirnov
statistics \cite{Abb}. One of  the one-sided statistics has the form
$$H^{+}_n= \sup_t \left[\Delta F_n(t)-\int_{-\infty}^{\infty} \left( \Delta
F_n(y) \right) dF_n(y)\right], $$ the second statistic $H_n^-$ is
introduced in a similar manner, and $H_n = \max(H_n^+, H_n^-).$

Formally, these statistics do not belong to the class of statistics
$SU_n^{\pm}$ and $SU_n,$  however the difference between them is
nonsignificant. Without loss of generality  we can assume that the
distribution of $X_i$ is uniform on $[-1,1]$. For any $t \in [-1,1]$
consider the family of statistics
$$H^{+}_n(t)=F_n(t)+F_n(-t)-\int_{-1}^1 (F_n(t)+F_n(-t)) dF_n(t)= n^{-2} \sum_{i,j=1}^{n}
\Psi(X_i,X_j;t)+n^{-1},$$ where the family of kernels $\Psi(\cdot \
; t)$  has the form
$$\aligned \Psi(x_1,x_2; t)&=\frac{1}{2} \Big (
I\{x_1 <t\}+I\{x_1<-t\}+ I\{x_2<t\}+I\{x_2<-t\}\\
&-2\cdot I\{x_1+x_2<0\}-1 \Big ).
\endaligned $$

 Now we introduce the auxiliary family of statistics
    $${\cal {H} }_n^{+}(t)={n \choose 2}^{-1}\sum_{1 \le i <j \le n}
\Psi(X_i,X_j; t).$$ Note that for any $t$ and any $n>1$
$$|H_n^{+}(t)-{\cal H}_n^{+}(t)| \le 7n^{-1},$$
so  the large deviation asymptotics for $\sup_t H_n^+(t)$ and
$\sup_t {\cal H}_n^+(t)$  are the same. It is easy to check that
$$\psi(z; t) \equiv E\left(\Psi(X_1,X_2; t)| X_1=z\right) =\frac{1}{2}
\left(I\{z<t\}+ I\{z<-t\}+t-1 \right),$$ and hence
$$ \sigma_{\psi}^2(t) = \frac14 ( t^2 - |t| +\frac 13), \ -1 \leq t
\leq 1. $$ Clearly $\psi_0^2 =\sup_{-1\leq t \leq 1} \sigma^2_{\psi}
(t) = 1/12.$

Let us turn to the upper bound and consider the uniform partition of
$[-1,1]$  into $2N$ parts using the nodes $t_k = k/N,\ k =-N,...,N.$
Obviously,
$$
\sup_{t_k \leq t < t_{k+1} } {\cal H }_n^+(t) \leq {\cal H
}_n^+(t_{k+1})+ F_n(t_{k+1}) - F_n(t_{k}).
$$
 By Lemma 5.1 this ensures  the monotonicity condition.
Therefore for some continuous function $v_4$  as $a \to 0$
$$\lim_{n \to \infty} n^{-1}\ln {\mathbb P} (H^{+}_n > a) =  v_4(a) = -\frac32
a^2 + O(a^3). $$  Similar results hold for $H_n^{-}$ and $H_n$.

5. Another test of symmetry is based on the characterization
established by Baringhaus and Henze \cite{BarHen}: {\it The common
distribution of two independent rv's $X$ and $Y$ is symmetric with
respect to zero iff $|X|$ and $|\max(X,Y)|$ have the same
distribution.}

Consider two edf's based on the sample $X_1,...,X_n.$ Let
$$ L_n(x) = n^{-1}\sum_{j=1}^n{I(|X_j| \le x)}, x\ge 0 $$ and
$$ G_n(x) ={n\choose 2}^{-1}\sum_{1\le j<k\le n} I(|\max(X_j,X_k)|\le x), x\ge 0.
$$ Following  \cite{BarHen}, consider the statistic
$$ BH_n^+ = \sup_{x\ge 0} [L_n(x) - G_n(x)].
$$  The statistics $BH_n^-$ and  $BH_n$ are
defined analogously. All these statistics are distribu\-ti\-on-free,
and we may assume that the observations are uniformly distributed on
$[-1,1]$. The statistic $ BH_n^+ $ admits the representation $$
BH_n^+ = {n\choose 2}^{-1} \sup_{0<t<1} \sum_{1\le j<k \le n}
\Xi(X_j,X_k; t),  $$ where for any $t\in (0,1)$
$$ \Xi (X_j,X_k; t) = I\{|\max(X_j,X_k)| \le t\} -
\frac{1}{2} (I\{|X_j| \le t\} + I\{|X_k| \le t\}).
$$

 Formally this family does not fit our theory, however
replacing edf $F_n$ by edf $L_n$ leads to minimal changes in the
proofs. Simple calculations show that for any $z \in [-1,1]$ and
$0<t<1$ the projections of the kernels have the form
$$ \xi(z; t) = E\left(\Xi(X_1,X_2; t)|X_1 =z \right) =
 \left\{\begin{array}{cc}t/2, & {-1 \le z<-t} \\
0, & {-t \le z \le t} \\ -t/2, & {t<z \le 1.}
\end{array} \right.$$
Consequently the variance function is given by
$$ \sigma^2_{\xi} (t) = \frac{t^2}{8} \int_{-1}^{1} ( I\{-1\leq x
<-t\} - I\{t < x \leq 1\})^2 dx = \frac14 t^2 ( 1-|t|), -1\leq t
\leq 1. $$ The maximum of this function is attained at $\pm \frac23$
and is equal to $\frac{1}{27}.$ Hence the large deviation
asymptotics has the form
$$ \lim_{n \to \infty} n^{-1} \ln {\mathbb P}(BH_n^+ \geq a) = v_5(a)=
-\frac{27}{8}a^2 + O(a^3), \ a  \to 0. $$ Similar statements hold
true for  the statistics $BH_n^-$ and $BH_n.$

6. Consider the famous characterization of normality due to G. Polya
\cite{pol}: {\it Let $X$ and $Y$ be i.i.d. rv's with zero mean. Then
$X $ and $(X+Y)/\sqrt{2}$ have the same distribution iff $X$ and $Y$
are normally distributed with some positive variance.}

The integral test of normality based on this characterization was
proposed in \cite{muli}. Let us construct the scale-invariant
Kolmogorov-type test comparing the usual edf and the $U$-empirical
df, based on $\frac{X_i + X_j}{\sqrt{2}}.$ We arrive at the
one-sided statistic
$$ PO^+_n = \sup_t \left({n \choose 2}^{-1} \sum_{1\leq i < j \leq n
} I\{X_i + X_j < t\sqrt{2}\} - n^{-1}\sum_{k=1}^n I\{X_k < t\}
\right ).
$$
 The statistics $PO_n^-$ and $PO_n$ are introduced similarly.
Our statistic $PO_n^{+}$ corresponds to the family of kernels
$$ \Psi(x_1,x_2; t) = I\{x_1+x_2 < t\sqrt{2} \} - \frac12( I\{x_1
<t\}+ I\{x_2 <t\}) $$ and satisfies the monotonicity condition.

Denote by $\cal {N}$ the df of the standard normal law. Then the
projections of the kernels are
$$ \psi(s; t) = E(\Psi(X_1,X_2; t)|X_2 =s)= {\cal {N}}(t\sqrt{2} -s) -
\frac12 {\cal {N}}(t) -\frac12 1\{s<t\}. $$ Consequently, the
variance function is $$  \sigma_{\psi}^2(t) =
\int_{-\infty}^{\infty} {\cal {N}}^2( t\sqrt{2} -s)d{\cal
{N}}(s)-\int_{-\infty}^{t}{\cal {N}}(t\sqrt{2} -s)d{\cal {N}}(s) +
\frac14 {\cal {N}}(t)- \frac14 {\cal {N}}^2(t) .
$$

The problem of finding the maximum of this function analytically is
difficult. However, its plot clearly shows that this maximum  is
attained at zero, and consequently is equal to 1/48.

\begin{center}
\includegraphics[width=3in,angle=270]{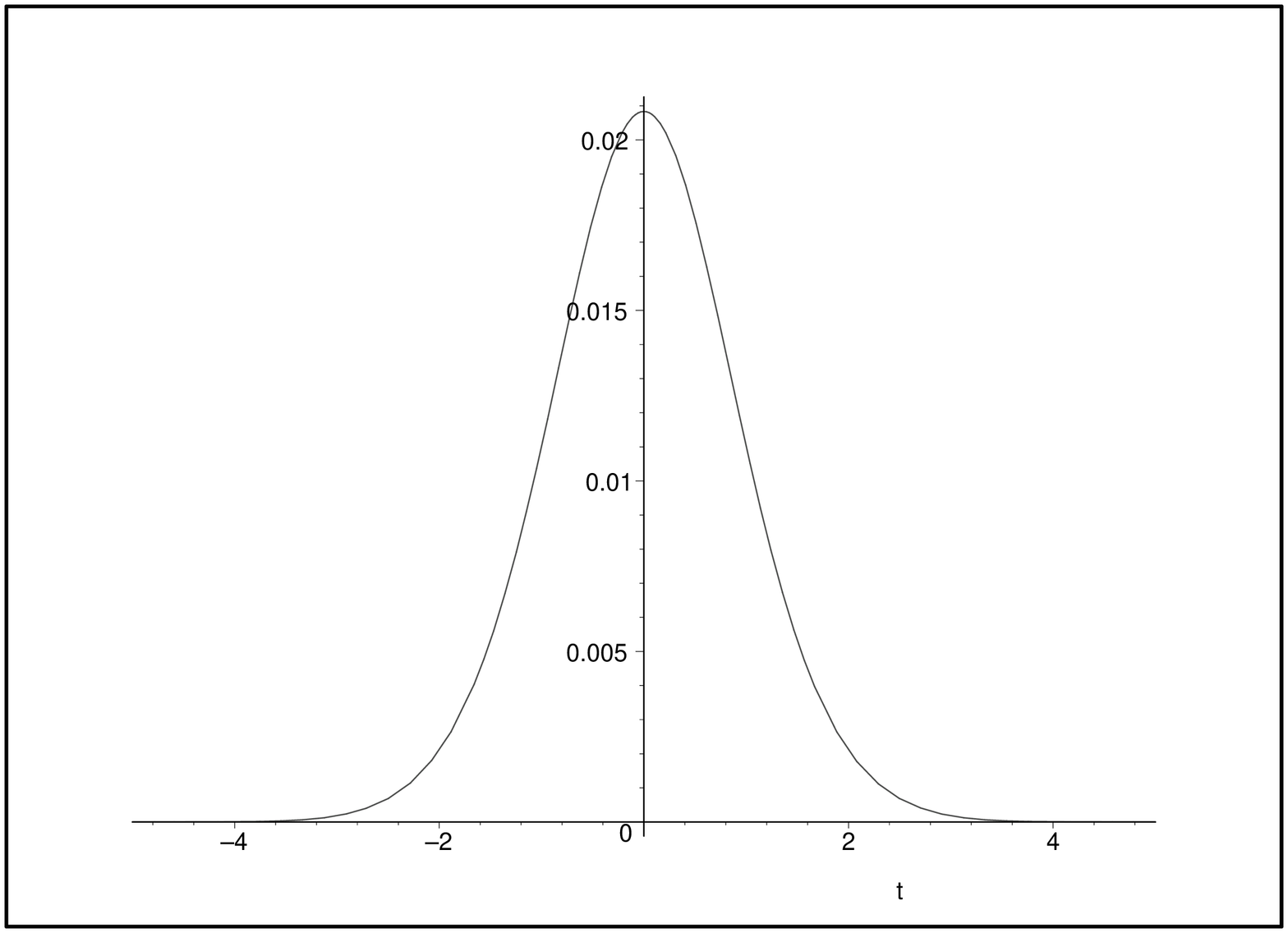}\\
\vspace{0.5cm} Fig.1 { \it Plot of the variance function for Polya
test.}
\end{center}

Therefore there exists a continuous function $v_6$ such  that
$$ \lim_{n\to \infty}n^{-1} \ln {\mathbb P}( PO_n > a) =  v_6(a) = -6 a^2
+ O(a^3),\ a \to 0.$$

\section{Calculation of local Bahadur efficiency}

The results on large deviations  allow us to calculate local Bahadur
efficiency of $U$-empirical Kolmogorov-Smirnov tests.

As an illustration we find the efficiency of Desu's test for some
parametric alternatives $F_{\theta}$ to the hypothesis of
exponentiality. The local exact slope, see \cite{Niki1995}, is the
main part as $\theta \to 0$  of the expression
$$4\sup_x |(1-F_{\theta}(x/2))^2 - (1- F_{\theta}(x))|^2.$$

Take, for example, the Weibull alternative with $1-F_{\theta}(x)=
\exp(- x^{1+\theta}), \theta \geq 0.$ Then the local exact slope is
equivalent to
$$
4\ln^2(2)\sup_x(x \exp(-x))^2 \ \theta^2 \approx 0.2601 \ \theta^2.
$$
The theoretical maximum (double Kullback-Leibler information, see
\cite{Bah1971}), is equal \cite{niki96}  to $\pi^2/6\ \cdot \theta^2
\approx 1.6449 \ \theta^2. $ Hence the local Bahadur efficiency of
Desu's test equals 0.1581. At the same time the efficiency of Desu's
test for Makeham alternative with the density $(1+\theta
(1-e^{-x}))\exp{(-x\!-\!\theta[x\!-\!(1\!-\!e^{\!-\!x})])}, x\geq 0,
\theta\geq 0, $ is much larger and is equal to 0.4938.

Other calculations of local Bahadur efficiency can be found in
\cite{Abb}, \cite{BH} and \cite{niki96}. It turns out that in some
cases $U$-empirical Kolmogorov tests have high efficiency and
perform well compared to some other goodness-of-fit tests. For
instance, the local efficiency of the sequence of statistics $H_n$
for testing symmetry under the normal shift alternative is equal to
$0.955.$ In the same problem the local efficiency of the sequence of
statistics $BH_n$ is equal to 0.75.

It would be interesting to construct new $U$-empirical tests of
Kolmogorov-Smirnov type and to calculate their efficiencies using
the large deviation results obtained  in this paper. We hope to
return to this question later.

\section*{Acknowledgements}

The author is thankful to Prof. A.I. Nazarov for reading the paper
and  making useful comments and to Prof. R. Serfling for valuable
advice.

\end{document}